\documentclass[lettersize,journal]{IEEEtran}
\usepackage{amsmath,amsfonts}
\usepackage{algorithmic}
\usepackage{algorithm}
\usepackage{array}
\usepackage[caption=false,font=normalsize,labelfont=sf,textfont=sf]{subfig}
\usepackage{textcomp}
\usepackage{stfloats}
\usepackage{url}
\usepackage{verbatim}
\usepackage{graphicx}
\usepackage{cite}
\usepackage{mathtools}
\usepackage{xcolor}
\usepackage{mathptmx}
\usepackage{amsmath}
\usepackage{amssymb}
\usepackage{amsthm}
\usepackage{float}
\usepackage{lipsum}
\usepackage{graphicx}
\DeclareMathOperator{\tr}{tr}
\renewcommand{\qedsymbol}{$\hfill\blacksquare$}

\newtheorem{thm}{Theorem}[section]

\newtheorem{prop}[thm]{Proposition}

\newtheorem{defn}{Definition}[section]

\newtheorem{rem}{Remark}
\newtheorem{assump}{Assumption}

\newcommand\blfootnote[1]{%
  \begingroup
  \renewcommand\thefootnote{}\footnote{#1}%
  \addtocounter{footnote}{-1}%
  \endgroup
}

\begin{document}

\textbf{This work has been submitted to the IEEE for possible publication. Copyright may be transferred without notice, after which this version may no longer be accessible.}
\title{Inverse linear-quadratic nonzero-sum differential games}

\author{Emin Martirosyan, Ming Cao\\
Engineering and Technology Institute Groningen\\ University of Groningen, Nijenborgh 4, Groningen, 9712CP, Netherlands}

\markboth{Journal of \LaTeX\ Class Files,~Vol.~14, No.~8, August~2021}%
{Shell \MakeLowercase{\textit{et al.}}: A Sample Article Using IEEEtran.cls for IEEE Journals}

\maketitle
\blfootnote{\textit{\textbf{This is a preprint of work that has been submitted to the IEEE for possible publication. Copyright may be transferred without notice, after which this version may no longer be accessible..}}}

\begin{abstract}
This paper addresses the inverse problem for Linear-Quadratic (LQ) nonzero-sum $N$-player differential games, where the goal is to learn parameters of an unknown cost function for the game, called observed, given the demonstrated trajectories that are known to be generated by stationary linear feedback Nash equilibrium laws. Towards this end, using the demonstrated data, a synthesized game needs to be constructed, which is required to be equivalent to the observed game in the sense that the trajectories generated by the equilibrium feedback laws of the $N$ players in the synthesized game are the same as those demonstrated trajectories. We show a model-based algorithm that can accomplish this task using the given trajectories. We then extend this model-based algorithm to a model-free setting to solve the same problem in the case when the system's matrices are unknown. The algorithms combine both inverse optimal control and reinforcement learning methods making extensive use of gradient descent optimization for the latter. The analysis of the algorithm focuses on the proof of its convergence and stability. To further illustrate possible solution characterization, we show how to generate an infinite number of equivalent games, not requiring to run repeatedly the complete algorithm. Simulation results validate the effectiveness of the proposed algorithms.
\end{abstract}

\begin{IEEEkeywords}
Inverse Differential Game, Inverse Optimal Control, Integral Reinforcement Learning, Continuous-time linear systems 
\end{IEEEkeywords}

\section{Introduction}

\IEEEPARstart{D}{ynamic} Game Theory is a branch of game theory that focuses on games where the strategies of the players can change over time \cite{basar1986dynamic}. Four features arise: the possible presence of multiple players (the number of players $N\geq2$), players' optimizing behavior, enduring consequences of decisions, and robustness against the changing environment \cite{engwerda_lq_nodate}. This dynamic aspect has gained significant attention in recent years, as many real-world problems modeled by games involve situations where the parameters of the game are constantly evolving \cite{basar1999dynamic}. For example, dynamic games can be used to model the competition of firms in a market, the evolution of political powers, and the interactions between populations in an ecosystem \cite{mailath06}. The typical dynamic games include differential games \cite{mailath06}, repeated games, and evolutionary games \cite{maynard1982evolutionary}. The study of these games has far-reaching implications in a range of fields such as economics \cite{sann08}, political science \cite{LEONG2010552}, engineering \cite{flad17}, \cite{mylv17}, \cite{gu08}, and biology \cite{now04}. Most of the literature has focused on determining the outcome of a game given the players' objective functions. Recently, interest grows in the inverse problem, where, given a player's game-playing behavior, one wants to reverse engineer the objective of this player.  

Inverse problems are particularly relevant  in guiding a game-playing system to some desired behavior. Inverse Reinforcement Learning (IRL), first introduced in \cite{Ng00algorithmsfor}, solves the inverse problem in a Markov Decision Process (MDP), using, e.g., maximum entropy methods \cite{ziebart_maximum_nodate}, \cite{ho2016generative}. Inverse Optimal Control (IOC), a closely related field with a long history, has focused on developing mathematical models and algorithms for inferring the objectives and constraints of a system in view of observed behavior. One of the earliest works in this area is the classic paper by Anderson in 1966 \cite{anderson1966inverse}, which introduced a linear-quadratic framework for inverse optimal control of linear systems. Further development of this framework leads to more results on IOC, e.g., \cite{menn18}, \cite{jean_inverse_2018}. IRL and IOC are concerned with similar problems, but differ in structure - the IOC aims to reconstruct an objective function given the state/action samples assuming dealing with a stable control system, while the IRL recovers an objective function using expert demonstration assuming that the expert behavior is optimal \cite{ab_azar_inverse_2020}.

Non-cooperative differential games were first introduced in \cite{isaacs_differential_1965} for zero-sum games. In this work, we consider a particular type of differential game - the LQ nonzero-sum game. This type of game is closely related to Linear Quadratic Regulator (LQR) problem -- the dynamics of the system are described by ordinary differential equations and the cost function is quadratic.  Thus, those methods used for solving IOC problems can be exploited for solving inverse differential games \cite{molloy_inverse_2017}. There are various works dedicated to the inverse problem for non-cooperative linear-quadratic differential games. Some of them use purely IRL approaches \cite{inga_solution_2019}, \cite{kopf_inverse_2017}, while others are based on IOC \cite{molloy_inverse_2020}, \cite{lian_robust_2021}.  

Our work considers LQ $N$-player differential games with heterogeneous players whose the control input matrices and cost function parameters are different. The solution for the considered type of game, namely its Nash equilibrium, is found via solving the so-called Algebraic Riccati Equation (AREs)  \cite{vamvoudakis_online_nodate},\cite{basar1999dynamic}. We exploit the result of \cite{li_1995} to accomplish this task. Further, instead of seeking the cost functions that, together with the dynamics, generated the demonstrated behavior, we look for an equivalent cost function that, together with the given dynamics, synthesizes a game that shares the same feedback laws with the original game. This can be done via model-based and model-free algorithms presented in this paper. The model-free algorithm, as an extension of the model-based version, is developed relying on he ideas of \cite{modares__2015}, \cite{jiang_computational_2012} and \cite{vamvoudakis_non-zero_2015} and using integral RL \cite{vrabie_adaptive_2009}. The extended algorithm possesses the same analytical properties as the model-based one. After characterizing the solution, we show that using the heterogeneity of the players, the output of the algorithms can be further adjusted allowing to generate an infinite number of equivalent games by exploiting such an algorithm again (thus, low computational costs).

The paper is structured as follows. Section \ref{PFsec} provides preliminary results on LQ nonzero-sum $N$-player differential games and formulates the problem addressed in the paper. In section \ref{MBsec}, we describe each step of the model-based algorithm. Section \ref{ANsec} is dedicated to the analysis of the algorithm; we show its convergence and stability and explain how to adjust the output of the algorithm via solution characterization. In section \ref{MFsec} we provide the model-free extension of the algorithm and show the equivalence of analytical results. Sections \ref{SIMsec} and \ref{CONsec} provide simulation results and conclusion, respectively.

\textit{Notations}: For a matrix $P\in\mathbb{R}^{m\times n}$, $P^k$, $P^{(k)}$ denote $P$ to the power of $k$, and the matrix $P$ at the $k$-th iteration, respectively. In addition, $P>0$, $P\geq0$, $P\leq0$, and $P<0$, denote positive (semi-)definiteness, and negative (semi-)definiteness of the matrix $P$, respectively. The notations $\{P_i\}_{i=1}^N$ and $\{P_{ij}\}_{i,j=1}^N$ denotes the sets of matrices $P_1,\dots,P_N$ and $P_{11},\dots,P_{1N},P_{21},\dots,P_{NN}$, respectively. The notation $\tr P$ denotes the trace of the matrix $P$. $I_k$ is the $k\times k$ identity matrix. 

\section{Problem Formulation}\label{PFsec}

This section introduces linear-quadratic (LQ) nonzero-sum differential games. We define stationary linear feedback Nash equilibrium (further referred to as NE). We clarify what an optimal behavior for the game is and introduce the inverse differential games.

\subsection{LQ Nonzero-sum Differential Game}

Consider a differential game with $N$ players, labeled by $1, ..., N$, under the  continuous time dynamics 
\begin{align}\label{stabdyn}
\begin{split}
    &\Dot{x}(t) = A x(t) + \sum_{i=1}^N B_i u_i (t), \quad i=1,\dots,N,\\
    &x(0) = x_0
\end{split}
\end{align}
where $x\in\mathbb{R}^n$ is the state and $u_i\in\mathbb{R}^{m_i}$ is the control input of players $i$; the plant matrix $A$, control input matrices $B_i$ have appropriate dimensions.

We consider that the players select their control to be
\begin{equation}
	u_ i(t) = F_i x(t),\quad i=1,\dots,N 
\end{equation}
where $F_i$ is an $m_i\times n$ time-invariant feedback matrix of player $i$. Further, to ease notations, we use $x(t)=x$, $u_i(t)=u_i$ for $i=1,\dots,N$.

We use $u_{-i} = (u_1,\dots,u_{i-1},u_{i+1},\dots,u_N)$ to denote an action profile of all the players except for player $i$. Within the game, player $i$ aims to find a controller $u_i$ that minimizes its cost function $J_i(x_0,u_i,u_{-i})$, which takes the quadratic form 
\begin{equation}\label{cost}
	J_i(x_0,u_i,u_{-i}) = \int_{0}^{\infty} \Big(x^\top Q_i x + \sum_{j=1}^N u_j^\top R_{ij} u_j\Big)\, dt,
\end{equation}
where $Q_i\in\mathbb{R}^{n\times n}$, $R_{ij}\in\mathbb{R}^{m_j\times m_j}$ are symmetric and $R_{ii}>0$ for $i,j=1,\dots,N$.

A Nash equilibrium $(u_i^*, u_{-i}^*)$ of the game is characterized by 
\begin{equation}
J_i(x_0,u_i^*, u_{-i}^*) \leq J_i(x_0,u_i, u_{-i}^*),\quad i=1,\dots,N.
\end{equation}
According to \cite[Theorem~8.5]{engwerda_lq_nodate}, for each player $i$ the cost function under the NE control inputs satisfies
\begin{equation}
J_i(x_0,u_i^*, u_{-i}^*)  = x_0^\top K_i x_0
\end{equation}
where $K_i$ is a symmetric matrix, sometimes referred to as the value matrix, satisfying the following Algebraic Riccati Equations (AREs)
\begin{align}\label{AREs}
    \begin{split}
        &A^\top K_i + K_i A + Q_i + \sum_{j=1}^N F_j^{\top} R_{ij} F_j -\\
        & \big(\sum_{j=1}^N F_j^{\top} B_j^\top \big) K_i - K_i \big(\sum_{j=1}^N B_j F_j\big) = 0,
    \end{split}
\end{align}
where $F_i$, for each $i=1,\dots,N$, is given by
\begin{equation}\label{eqfdblw}
    F_i = R_{ii}^{-1} B_i^\top K_i,
\end{equation}
and the control trajectories are 
\begin{equation}
u_i^* = -F_i x = - R_{ii}^{-1} B_i^\top K_i x,\quad i=1,\dots,N.
\end{equation}

We restrict the set of admissible controller matrices $(F_1,\dots,F_N)$ to the following set
\begin{equation}\label{stab}
\mathcal{F} = \{(F_1,\dots,F_N)| A + \sum_{j=1}^N B_j F_j\quad\text{is stable}\},
\end{equation}
since $(u_1^*,\dots,u_N^*)$ need to stabilize trajectories to qualify as the NE equilibrium in this game \cite{engwerda_lq_nodate}. This restriction is essential because, as shown in \cite{mage76}, without this restriction it is possible to construct an example where a non-stabilizing feedback yields a lower cost for one of the player while other players stick to the stabilizing feedback law. Thus, besides satisfying \eqref{AREs}, $K_i$ for $i=1,\dots,N$ should also be stabilizing to lead to an NE \cite{engwerda_lq_nodate}. Thus, the system \eqref{stabdyn} in the LQ differential games is always assumed to be stabilizable, i.e., $(A,[B_1,\dots,B_N])$ is stabilizable. 

\subsection{Inverse LQ Nonzero-sum Differential Game}

We formulate the inverse problem for LQ nonzero-sum differential games in this subsection. 

Consider an LQ differential game (referred to as the observed LQ game) with continuous-time system dynamics
\begin{equation}\label{demdyn}
    \dot{x}_d = A x_d + \sum_{i=1}^N B_i u_{i,d},\quad x_d(0) = x_{0,d}
\end{equation}
where $x_d\in\mathbb{R}^n$, $u_{i,d}\in\mathbb{R}^{m_i}$ are the demonstrated NE trajectories of the observed LQ game with $u_{i,d}$ being the trajectory of player $i$ for $i=1,\dots,N$; $A$, $B_i$ have appropriate dimensions. The cost functions of the game have the following known quadratic structure
\begin{equation}\label{cost1}
	J_i(x_0, u_i, u_{-i}) = \int_{0}^{\infty} \Big(x^\top Q_{i,d} x + \sum_{j=1}^N u_j^\top R_{ij,d} u_j)\Big) dt,
\end{equation}
with the \emph{unknown} symmetric matrices $Q_{i,d}$ and $R_{ij,d}$ where $R_{ii,d}>0$ for $i,j=1,\dots,N$. 
Considering that $(x_d, u_{1,d}, \dots,u_{N,d})$ are NE trajectories, we have
\begin{equation}\label{trgtfdbkl}
    u_{i,d} = - F_{i,d} x_d = - R_{ii,d}^{-1} B_i^\top K_{i,d} x_d,
\end{equation}
where $K_{i,d}$ is the stabilizing symmetric solution of the following AREs
\begin{align}\label{invpAREs}
\begin{split}
    &A^\top K_{i,d} + K_{i,d} A + Q_{i,d} + \sum_{j=1}^N F_{j,d}^{\top} R_{ij,d} F_{j,d} -\\
    & \big(\sum_{j=1}^N F_{j,d}^{\top} B_j^\top \big) K_{i,d} - K_{i,d} \big(\sum_{j=1}^N B_j F_{j,d}\big) = 0.
    \end{split}
\end{align}
\begin{rem}
    Note that we do not make any assumption on stabilizability of the system because it follows from the existence of demonstrated NE trajectories.
\end{rem}

We use the $(A,\{B_i\}_{i=1}^N,\{Q_{i,d}\}_{i=1}^N,\{R_{ij,d}\}_{i,j=1}^N)$ tuple to describe an LQ differential game with the dynamics' matrices $A,B_1,\dots,B_N$ and the cost function parameters $Q_{1,d},\dots,Q_{N,d}$ and $R_{11,d},\dots,R_{1N,d},R_{21,d},\dots,R_{NN,d}$. 

\begin{defn}{(Equivalent Game).}\label{defeq}
The $(A,\{B_i\}_{i=1}^N,\{Q_i\}_{i=1}^N,\{R_{ij}\}_{i,j=1}^N)$ game is said to be equivalent to the observed game $(A,\{B_i\}_{i=1}^N,\{Q_{i,d}\}_{i=1}^N,\{R_{ij,d}\}_{i,j=1}^N)$ if its AREs \eqref{AREs} has a stabilizing solution $\{K_i\}_{i=1}^N$ such that $R_{ii}^{-1} B_i^\top K_i = R_{ii,d}^{-1} B_i^\top K_{i,d}$ (i.e., $F_i= F_{i,d}$) where $\{K_{i,d}\}_{i=1}^N$ is a solution of AREs \eqref{invpAREs} associated with the observed game.
\end{defn}
\noindent In other words, the games are equivalent if they share the same equilibrium feedback laws $F_{i,d} = F_i$ for all player $i=1,\dots,N$.

Now, we are ready to formulate the inverse problem to be addressed in this paper.\\
\textbf{Inverse Differential Game Problem}: Given the demonstrated trajectories $(x_d,u_{1,d},\dots,u_{N,D})$ of the observed game $(A,\{B_i\}_{i=1}^N,\{Q_{i,d}\}_{i=1}^N,\{R_{ij,d}\}_{i,j=1}^N)$, find the cost function parameters $\{Q_i,R_{ij}\}_{i,j=1}^N$ that synthesize a game $(A,\{B_i\}_{i=1}^N,\{Q_i\}_{i=1}^N,\{R_ij\}_{i,j=1}^N)$ which is equivalent to the observed game.

We solve the problem using model-based and model-free algorithms presented in the following sections.

\section{Model-based Inverse Reinforcement Learning Algorithm}\label{MBsec}

This section describes the algorithm that uses the demonstrated equilibrium trajectories $(x_d,\{u_{i,d}\}_{i=1}^N)$ generated by the \emph{known} dynamics $(A,\{B_i\}_{i=1}^N)$ for learning a set of cost function parameters equivalent to $(Q_{i,d},R_{ij,d})$ for $i,j=1,\dots,N$.

The algorithm consists of the following steps - firstly, we use the demonstrated data to estimate the set of target feedback laws $\hat{F}_i = F_{i,d}$ that are supposed to be a set of equilibrium feedback laws both for the original game and the one generated by the algorithm. The next step is the initialization of the parameters $(Q_i^{0)}$, $R_{ij})$ for $i,j=1,\dots,N$. Note that the algorithm only updates $Q_i$'s parameters while $R_{ij}$'s remain the same during the iterative procedure. Then, using the initialized parameters and the known dynamics, we calculate the set of stabilizing solutions of the resulting set of AREs and the corresponding feedback laws. Using the initialized feedback laws $F^{(k)}_i$, we apply the gradient descent method \cite{Bertsekas/99} to update $K_i^{(k)}$ in the direction of the minimization of the difference between $F_i^{(k)}$ and $\hat{F}_i$ for $k=0,1,\dots$. After each iteration, using the inverse optimal control \cite{haddad08}, we update $Q_i^{k+1}$ substituting the result of the gradient descent update $K_i^{(k+1)}$. 

The model-based algorithm requires to know matrices of the game dynamics. Hence, in this section we make the following assumption.
\begin{assump}
The game dynamics matrices $(A,B_1,\dots,B_N)$ are known.
\end{assump}
\subsection{Feedback Law Estimation}

In this step we aim to track the difference between the the current iteration $k$ feedback laws and the desired ones. Using the observed data $(x_d,\{u_{i,d}\}_{i=1}^N)$, we derive the estimation $\hat{F}_i$ of the target feedback law $F_{i,d}$ by applying the batch least-square (LS) method \cite{devore2015probability}. To implement the estimation procedure we sample the demonstrated trajectories to obtain 
\begin{align}
\begin{split}
    \hat{x}_d &= [x_d(t_1),\dots,x_d(t_s)]\in\mathbb{R}^{n\times s},\\
    \hat{u}_{i,d} &= [u_{i,d}(t_1),\dots,u_{i,d}(t_s)]\in\mathbb{R}^{m_i\times s},
\end{split}
\end{align}
for $i = 1,\dots,N$ where $s\geq n$, $s\in\mathbb{Z}_+$. Using \eqref{trgtfdbkl}, we estimate $\hat{F}_i$ by calculating 
\begin{equation}\label{LQest}
    \hat{F}_i = - \hat{u}_{i,d} \hat{x}_d^\top (\hat{x}_d \hat{x}_d^\top)^{-1},
\end{equation}
for $i=1,\dots,N$. Note that the sampling should guarantee that $\hat{x}_d \hat{x}_d$ is full rank, i.e., the rank should be $n$.
\subsection{Initialized Game}

In the next step, we generate an initial set of parameters $\{Q^{(0,0)}_i, \{R_{ij}\}_{j=1}^N\}_{i=1}^N$. Together with the matrices $(A,\{B_i\}_{i=1}^N)$ we have a nonzero-sum linear quadratic differential game. To find the equilibrium set $\{F_i^{(\infty,0)}\}_{i=1}^N$ for the generated game, one needs to solve the following set of equations
\begin{align}
\begin{split}
    &A^\top K_i^{(0,0)} + K_i^{(0,0)} A + Q_i + \sum_{j=1}^N F_j^{(0,0)\top} R_{ij} F_j^{(0,0)} -\\
    & \big(\sum_{j=1}^N F_j^{(0,0)\top} B_j^\top \big) K_i^{(0,0)} - K_i^{(0,0)} \big(\sum_{j=1}^N B_j F_j^{(0,0)}\big) = 0,
\end{split}
\end{align}
where $F_i^{(0,0)} =  R_{ii}^{-1} B_i^\top K_i^{(0,0)}$. This set of AREs can be solved using a modified version, for the multiplayer case, of the algorithm of the Lyapunov Iterations presented in \cite{li_1995}. The algorithm includes initialization of $F_i^{(0,0)}$ that should form stable dynamics, i.e., 
\begin{equation}
    A - \sum_{i=1}^N B_i F_i^{(0,0)}\quad\text{is stable}.
\end{equation}
However, since $\{\hat{F}_i\}_{i=1}^N$ is derived using the estimation procedure and known to be a set of equilibrium feedback laws, one can skip the initialization step for solving the set of AREs by setting $F_i^{(0,0)} = \hat{F}_i$. Thus, the algorithm used to solve initialized game is the following 
\begin{align}\label{lyapiter}
\begin{split}
    &(A - \sum_{j=1}^N B_j F_j^{(k,0)})^\top K_i^{(k+1,0)} + K_i^{(k+1,0)} (A - \sum_{j=1}^N B_j F_j^{(k,0)}) =\\
    & -\Tilde{Q}_i = -(Q_i + \sum_{j=1}^N F_j^{(k,0)\top} R_{ij} F_j^{(k,0)}),\\
    &F_i^{(k+1,0)} = R_{ii}^{-1} B_i^\top K_i^{(k+1,0)}
\end{split}
\end{align}
for iterations $k=0,1,\dots$. The procedure continues until $\lVert K_i^{(k+1,0)} - K_i^{(k,0)}\rVert \leq \epsilon_i$ where $\epsilon_i$ is some positive constant for $i=1,\dots, N$. From \cite{li_1995}, we know that under some mild conditions, the algorithm converges to a set of \emph{positive definite stabilizing} solutions $\{K_i^{(\infty,0)}\}_{i=1}^N$, i.e., 
\begin{equation}\label{solinit}
    \lim_{k\rightarrow\infty} K_i^{(k,0)} = K_i^\infty,
\end{equation}
where $K_i^{(\infty,0)}$ are such that $A - \sum_{i=1}^N K_i^{(\infty,0)}$ is stable. Because \eqref{lyapiter} is a set of the Lyapunov Equations, the conditions that ensure the uniqueness of the set are the following 
\begin{equation}
    Q_i >0,\, R_{ii} > 0,\, R_{ij} \geq 0
\end{equation}
for $i,j = 1,\dots,N$ and $i\neq j$.

Thus, we set the initialized parameters $Q_i^{(0)}$ and $R_{ij}$ as positive definite for $i=j$ and positive semi-definite for $i\neq j$, $i,j=1,\dots,N$. Note that further, in Section \ref{solcharact}, dedicated to the solution characterization, we show that $R_{ij}$ and the resulting $Q_i$'s can be adjusted relaxing the imposed constraint.  

After solving the initialized game, we set $K_i^{(\infty,0)} = K_i^{(0)}$ and correspondingly $F_i^{(\infty,0)} = F_i^{(0)}$.

\subsection{Gradient Descent Update}
 In this section we present the way we track the difference between the estimated controller $\hat{F}_i$ and $F_i^{(p)} = R_{ii}^{-1} B_i^\top K_i^{(p)}$, where $p=0,1,2,\dots$ is an iteration step and $K_i^{(0)}$ are the solution of the initialized problem \eqref{solinit} for $i=1,\dots,N$. This step is performed using the gradient descent algorithm \cite{Bertsekas/99}. We define the following functions
\begin{equation}\label{diff}
    d_i^{(p)} (K_i) \coloneqq F_i^{(p)} - \hat{F}_i = R_{ii}^{-1} B_i^\top K_i^{(p)} - \hat{F}_i
\end{equation}
which track the difference between the target feedback law $\hat{F}_i$ and the current iteration feedback law $F_i^{(p)}$ for player $i=1,\dots,N$. Next, we introduce the function $D_i^{(p)}$ of $K_i$ as follows
\begin{equation}\label{trdiff}
    D_i^{(p)} (K_i) = \tr{\big(d_i^{(p)\top} d_i^{(p)}\big)}\geq0,\quad i=1,\dots,N
\end{equation}
which we minimize with respect to $K_i^{(p)}$. The update rule is the following 
\begin{equation}\label{updatemb}
    K_i^{(p+1)} = K_i^{(p)} - \alpha_i \frac{\partial D_i^{(p)}}{\partial K_i},
\end{equation}
for $i = 1,\dots, N$ where $\alpha_i\geq 0$ is the learning rate for player $i$. Considering \eqref{eqfdblw}, \eqref{diff} and \eqref{trdiff}, we compute the partial derivative as follows
\begin{align}
\begin{split}
    \frac{\partial D_i^{(p)}}{\partial K_i} &= K_i^{(p)} B_i R_{ii}^{-1}  R_{ii}^{-1} B_i^\top + B_i R_{ii}^{-1} R_{ii}^{-1} B_i^\top K_i^{(p)} - \\
    & \hat{F}_i^\top R_{ii}^{-1} B_i^\top - B_i R_{ii}^{-1} \hat{F}_i\\
    & = (F_i^{(p)} - \hat{F}_i)^\top R_{ii}^{-1} B_i^\top + B_i R_{ii}^{-1} (F_i^{(p)} - \hat{F}_i)\\
    & = d_i^{(p)\top} R_{ii}^{-1} B_i^\top + B_i R_{ii}^{-1} d_i^{(p)}.
\end{split}
\end{align}
At each iteration $p=0,1,\dots$ we have bounded $d_i^{(p)}$ for $i=1,\dots,N$ because $d_i^{(0)} = F_i^{(\infty,0)} - \hat{F}_i$ where $F_i^{(0)}$ is the solution of the initialized problem and the following   
\begin{equation}
    C_i = \lVert d_i^{(0)} \rVert > \lVert d_i^{(1)} \rVert > \dots \geq 0
\end{equation}
as the result of the minimization procedure where $C_i \geq 0$ is a constant for $i=1,\dots,N$.

\subsection{Inverse Update of the Parameters}
After the update \eqref{updatemb}, we use $K_i^{(p+1)}$ to evaluate $Q_i^{(p+1)}$ for $i=1,\dots,N$. This is done via substituting the derived values into 
\begin{align}
\begin{split}
    & Q_i^{(p+1)} = -A^\top K_i^{(p+1)} - K_i^{(p+1)} A - \\
    & \sum_{j=1}^N F_j^{(p+1)\top} R_{ij} F_j^{(p+1)} + \big(\sum_{j=1}^N F_j^{(p+1)\top} B_j^\top \big) K_i^{(p+1)} + \\ 
    & K_i^{(p+1)} \big(\sum_{j=1}^N B_j F_j^{(p+1)}\big),
\end{split}
\end{align}
where $F_i^{(p+1)} = R_{ii}^{-1} B_i^\top K_i^{(p+1)}$ for $i=1,\dots,N$.

The described iterative procedure is repeated till for some $\delta_i$, $0 \leq D_i^{(p)} \leq \delta_i$ is achieved where $\delta_i$ are desired precision measures that describe how close the generated parameters are to the desired result for each player $i=1,\dots,N$. The resulting $Q_i^*$, together with the initialized $R_{ij}$ for $i,j=1,\dots,N$ and the known dynamics $(A,\{B_i\}_{i=1}^N)$ form an equivalent LQ nonzero-sum game as described in Definition \ref{defeq}. Hence, we have a new set of Algebraic Riccati Equations
\begin{align}
\begin{split}
    & Q_i^* = -A^\top K_i^* - K_i^* A - \sum_{j=1}^N F_j^{*\top} R_{ij} F_j^* + \\
    &\big(\sum_{j=1}^N F_j^{*\top} B_j^\top \big) K_i^* + K_i^* \big(\sum_{j=1}^N B_j F_j^*\big),
\end{split}
\end{align}
where $K_i^*$ is the final result of \eqref{updatemb} and $F_i^* = R_{ii}^{-1} B_i^\top K_i^* = \hat{F}_i$ for $i=1, \dots, N$.

\begin{rem}
    From the complexity point of view, the demanding parts of the algorithm are finding solutions of the game with the initialized parameters $\{Q_i^{(0)}, R_{ij}\}_{i,j=1}^N$ and matrix multiplication done in the following steps. Implementing the Lyapunov Iterations with respect to $K_i\in\mathbb{R}^{n\times n}$ usually has complexity $\mathcal{O}(n^3)$ \cite{golub2013matrix}. The steps of the algorithm that require performing matrix multiplication via standard methods have complexity $\mathcal{O}(n^3 + n^2 m + n m^2)$ where $m=\max(m_1,\dots,m_N)$. Hence, the overall computational complexity is $\mathcal{O}(n^3 + n^2m + n m^2)$.
\end{rem}

\begin{rem}\label{remstep}
In fact, the implementation of \textbf{Algorithm 1} does not necessarily require the iterative update of $Q_i^{(p+1)}$ in step 6. This update might be done only \emph{once} after the desired precision $\delta_i$ is achieved, i.e., after getting $K_i^{(p+1)}$ in step 5 such that $\tr(d_i^{(p)\top} d_i^{(p)}) < \delta_i$ for $i=1,\dots,N$. This would reduce the computational cost of the algorithm. On the other hand, steps $4$, $5$ and $6$ can be combined by substituting \eqref{diffcalc} and \eqref{UPD} into \eqref{IOC}.
\end{rem}

\begin{rem}
    Step $3$ is only necessary to derive solutions for the game with the initialized parameters. Suppose we are given a set of game parameters $\{Q_i^\prime,R_{ij}^\prime\}_{i,j=1}^N$ and the solution for that game $\{K_{i}^\prime,F_i^\prime\}_{i=1}^N$ is known to have the same dynamics as the observed game. Then, considering Remark \ref{remstep}, we only need to perform iterative optimization via steps 4-5 and a single update in step 6. The same applies if $A$ is known to be stable. In that case, one can skip Step $3$ and set $K_i^{(0)} = \mathbf{0}_{n\times n}\in\mathbb{R}^{n\times n}$ and $F_i^{(0)} = \mathbf{0}_{m_i\times n}\in\mathbb{R}^{m_i\times n}$ for $i=1,\dots,N$ where $\mathbf{0}$ denotes a zero matrix of particular dimension.
\end{rem}

\begin{algorithm}\label{alg1}
	\caption{Model-based Inverse Reinforcement Learning Algorithm}
	\begin{enumerate}
	    \item Initialize $R_{ii}> 0$, $R_{ij}\geq 0$ and $Q_i^{(0)}> 0$ for $i,j = 1,\dots,N$, $i\neq j$. Sample data from demonstrated $(x,\{u_{i,d}\}_{i=1}^N)$ to generate $(\hat{x},\{\hat{u}_{i,d}\}_{i=1}^N)$. Set $k=0$ and $p=0$.
	    \item Derive the estimation of  $F_{i,d}$ using the sampled data as
	    \begin{equation}
	        \hat{F}_i =  -\hat{u}_{i,d}\hat{x}_d^\top(\hat{x}_d \hat{x}_d^\top)^{-1}.
	    \end{equation}
	    \item Set $F_i^{(0,0)} = \hat{F}_i$ , compute  $K_i^{(k+1,0)}$ from 
    \begin{align}\label{mbinitialsol}
    \begin{split}
        &(A - \sum_{i=1}^N B_i F_i^{(k,0)})^\top K_i^{(k+1,0)} + K_i^{(k+1,0)} (A - \sum_{i=1}^N B_i F_i^{(k,0)}) =\\
        & -(Q_i^{(0)} + \sum_{j=1}^N F_j^{(k,0)\top} R_{ij} F_j^{(k,0)}),
    \end{split}
    \end{align}
    update 
    \begin{equation}\label{initialFL}
        F_i^{(k+1,0)} = R_{ii}^{-1} B_i^\top K_i^{(k+1,0)},
    \end{equation}
	    and set $k = k+1$ till $\lVert K_i^{(k+1,0)} - K_i^{(k,0)} \rVert < \epsilon_i$ where $\epsilon_i$ is a small positive constant for $i=1,\dots, N$. 
	    \item Set $K_i^{(0)} = K_i^{(k+1,0)}$, $F_i^{(0)} = F_i^{(k+1,0)}$. Evaluate the difference
	    \begin{equation}\label{diffcalc}
	    d^{(p)}_i = F_i^{(p)} - \hat{F}_i.
	    \end{equation}
	    \item Update $K^{(p+1)}_i$ and $F^{(p+1)}_i$ for $i=1,\dots,N$ as
	    \begin{align}\label{UPD}
     \begin{split}
	        &K_i^{(p+1)} = K_i^{(p)} - \alpha_i \Big(d^{(p)\top}_i R_{ii}^{-1} B_i^\top + B_i R_{ii}^{-1} d_i^{(p)}\Big),\\
            & F^{(p+1)}_i = R_{ii}^{-1} B_i^\top K_i^{(p+1)}. 
	 \end{split}   
     \end{align}
	    \item Perform evaluation of $Q^{(p+1)}$ as         
     \begin{align}\label{IOC}
     \begin{split}
        & Q_i^{(p+1)} = -A^\top K_i^{(p+1)} - K_i^{(p+1)} A - 
         \sum_{j=1}^N F_j^{(p+1)\top} R_{ij} F_j^{(p+1)} \\ 
        & + \big(\sum_{j=1}^N F_j^{(p+1)\top} B_j^\top \big) K_i^{(p+1)} + K_i^{(p+1)} \big(\sum_{j=1}^N B_j F_j^{(p+1)}\big).
    \end{split}
    \end{align}
	    \item Set $p = p + 1$. Perform steps $4$-$6$ till $\tr(d_i^{(p)\top} d_i^{(p)}) < \delta_i$ where $\delta_i$ is a small positive constant for $i=1,\dots,N$.
	\end{enumerate}
\end{algorithm}

\section{Analysis of the Model-based Algorithm}\label{ANsec}
This section is dedicated to the analysis of the model-based algorithm -- \textbf{Algorithm 1}. Firstly, we show the convergence of the algorithm. Next, we prove that the output of the algorithm to solve the problem, i.e., $F_{i,d}$, the feedback laws used to generate the equilibrium trajectories $(x_d,\{u_{i,d}\}_{i=1}^N$, are equilibrium trajectories for the synthesized game. In the end, we give the characterization of the solutions that allows to create other equivalent games. 

We need to introduce the following notations
\begin{align}
    A_{cl}^{(k,0)} &\coloneqq A - \sum_{i=1}^N B_i F_i^{(k,0)},\\
    g_i (d_i^{(p)}) &\coloneqq d^{(p)\top}_i R_{ii}^{-1} B_i^\top + B_i R_{ii}^{-1} d_i^{(p)} = g_i^{(p)},
\end{align}
where $g_i^{(p)}$ is the symmetric matrix for $p=0,1,\dots$ and $i=1,\dots,N$.

\subsection{Convergence Analysis}
The result on the convergence is formulated in the theorem below. 

\begin{thm}\label{convtheorem}
In \textbf{Algorithm 1}, the state reward parameters $Q_i^{(p)}$ converge to $Q_i^*$ for $i=1,\dots,N$. Furthermore, $Q_i^*$ together with the initialized $R_{ij}$, $i,j=1,\dots,N$, and the dynamics matrices $(A,\{B_i\}_{i=1}^N)$ form AREs with the stabilizing solution $K_i^*$ such that 
\begin{equation}\label{alggoal}
    R^{-1}_{ii} B_i^\top K_i^* = R^{-1}_{ii,d} B_i^\top K_{i,d} = F_{i,d}.
\end{equation}
\end{thm}
\noindent\textbf{Proof}. After the initialization procedure, we get $\{F_i^{(\infty,0)}\}_{i=1}^N$ such that $A_{cl}^{(\infty,0)}$ is stable. Consider the update rule \eqref{UPD}. The gradient descent update drives the initialized $F^{(0)}_i = F_i^{(\infty,0)}$ to the estimation of the target feedback law $\hat{F}_i$ for $i=1,\dots,N$. Hence, the function that is optimized satisfies 
\begin{equation}
    0\leq D_i^{(p+1)} < D_i^{(p)},\quad i=1,\dots,N,\,p=0,1,\dots. 
\end{equation}
Thus, the following can be deduced
\begin{align}
\begin{split}
    &\lim_{p\to\infty} D_i^{(p)} = 0,\quad \lim_{p\to\infty} d_i^{(p)} =0\\
    &\lim_{p\to\infty} g^{(p)}_i = 0,\quad i=1,\dots,N.
\end{split}
\end{align}
Thus, 
\begin{equation}
    \lim_{p\to\infty} K^{(p+1)}_i = \lim_{p\to\infty} (K^{(p)}_i -\alpha_i g_i^{(p)}) = \lim_{p\to\infty} K_i^{(p)}
\end{equation}
and, since $\hat{F}_i = F_{i,d}$, one can conclude 
\begin{equation}\label{stabsol}
     \lim_{p\to\infty} R^{-1}_{ii} B_i^\top K^{(p)}_i = \lim_{p\to\infty} F^{(p)}_i = F_{i,d} = R_{ii,d}^{-1} B_i^\top K_{i,d}.  
\end{equation}
for $i=1,\dots,N$.

The result of the convergence is denoted by $K_i^*$ for $i=1,\dots,N$. Substituting $F_i^{(p+1)} = R_{ii}^{-1} B_i^\top K_i^{(p+1)}$ and the gradient descent update \eqref{UPD} in the form
\begin{equation}
    K_i^{(p+1)} = K_i^{(p)} - \alpha_i g_i^{(p)}
\end{equation}
into \eqref{IOC}, we get 
\begin{align}
\begin{split}
            & Q_i^{(p+1)} = A^\top (K_i^{(p)} - \alpha_i g_i^{(p)}) + (K_i^{(p)} - \alpha_i g_i^{(p)}) A +\\
         & \sum_{j=1}^N (K_i^{(p)} - \alpha_i g_i^{(p)}) B_j R_{jj}^{-1} R_{ij} R_{jj}^{-1} B_j^\top (K_i^{(p)} - \alpha_i g_i^{(p)})-\\ 
        & \big(\sum_{j=1}^N (K_i^{(p)} - \alpha_i g_i^{(p)}) B_j R_{jj}^{-1} B_j^\top \big) (K_i^{(p)} - \alpha_i g_i^{(p)}) - \\
        & (K_i^{(p)} - \alpha_i g_i^{(p)}) \big(\sum_{j=1}^N B_j R_{jj}^{-1} B_j^\top (K_i^{(p)} - \alpha_i g_i^{(p)})\big).
\end{split}
\end{align}
Taking the limit and using $F_i^{(p)} = R_{ii}^{-1} B_i^\top K_i^{(p)}$, we get
\begin{align}
    \begin{split}
        & \lim_{p\to\infty} Q_i^{(p+1)} = \lim_{p\to\infty} ( A^\top K_i^{(p)} + K_i^{(p)} A + 
         \sum_{j=1}^N F_j^{(p)\top} R_{ij} F_j^{(p)} \\ 
        &- \big(\sum_{j=1}^N F_j^{(p+1)\top} B_j^\top \big) K_i^{(p+1)} - K_i^{(p)} \big(\sum_{j=1}^N B_j F_j^{(p}\big)).
    \end{split}
\end{align}
and, as a result, 
\begin{equation}
    \lim_{p\to\infty} Q_i^{(p+1)} = \lim_{p\to\infty} Q_i^{(p)},\quad i=1,\dots,N.
\end{equation}
Denoting the result of convergence as $Q_i^*$, we obtain
\begin{align}\label{AREsConv}
    \begin{split}
               & Q_i^* = A^\top K_i^* + K_i^* A + 
         \sum_{j=1}^N F_j^{*\top} R_{ij} F_j^* \\ 
        &- \big(\sum_{j=1}^N F_j^{*\top} B_j^\top \big) K_i^* - K_i^* \big(\sum_{j=1}^N B_j F_j^*\big),
    \end{split}
\end{align}
for $i=1,\dots,N$. Thus, we conclude that $\{K^*_i\}_{i=1}^N$ is the solution set for the AREs associated with $\{Q_i^*,R_{ij}\}_{i,j=1}^N$ where $R_{ij}$ are initialized parameters in \textit{Step 1}. Moreover, from \eqref{stabsol}, once concludes that it is a stabilizing solution set.\qedsymbol

\subsection{Stability Analysis}
In this section, we show that the output of the algorithm is an equivalent game to the game that has the demonstrated NE trajectories, i.e., $(x_d,u_{i,d})$ for $i=1,\dots,N$.

Firstly, we need to present the following result, extended for the multiplayer case on LQ nonzero-sum differential games from \cite{engwerda_lq_nodate}. 
\begin{thm}\label{ENGthr}
    Let $(K_1,\dots,K_N)$ be a symmetric stabilizing solution of equations \eqref{AREs} and define $F_i^* \coloneqq R_{ii}^{-1} B_i^\top K_i$ for $i=1,\dots,N$. Then $(F_1^*,\dots,F_N^*)$ is the feedback NE. Conversely, if $(F_1^*,\dots,F_N^*)$ is the feedback NE, there exists a symmetric stabilizing solution $(K_1,\dots,K_N)$ of equations \eqref{AREs} such that $F_i^* = R_{ii}^{-1} B_i^\top K_i$ for $i=1,\dots,N$.
\end{thm}
Finally, one can conclude the following for the proposed algorithm. 
\begin{thm}
    Given the demonstrated trajectories $(x_d, u_{i,d})$ for $i=1,\dots,N$ generated by a game $(A,\{B_i\}_{i=1}^N,\{Q_{i,d}\}_{i=1}^N,\{R_{ij,d}\}_{i,j=1}^N)$ described in Section \ref{PFsec}, the output of \textbf{Algorithm 1} is the tuple $(\{Q_i^*\}_{i=1}^N,\{R_{ij}\}_{i,j=1}^N)$ which combined with the known dynamics matrices $(A,\{B_i\}_{i=1}^N)$, synthesizes a game equivalent to $(A,\{B_i\}_{i=1}^N,\{Q_{i,d}\}_{i=1}^N,\{R_{ij,d}\}_{i,j=1}^N)$, i.e., $F^*_i = F_{i,d}$ for $i=1,\dots,N$.
\end{thm}
\noindent\textbf{Proof}. From \eqref{AREsConv} we know that $K^*_i$, $i=1,\dots,N$ is the solution for AREs with the parameters $(\{Q_i^*\}_{i=1}^N,\{R_{ij}\}_{i,j=1}^N)$ and dynamics $(A,\{B_i\}_{i=1}^N)$. From Theorem \ref{convtheorem}, we know that $F_{i,d} = R_{ii}^{-1} B_i^\top K_i^* = F_i^*$. Since $\{F_{i,d}\}_{i=1}^N$ is the set of stabilizing feedback laws, $K_*$ is the set of stabilizing solutions for AREs with parameters generated by the algorithm and, as a result of Theorem \ref{ENGthr}, one conclude that $\{F_i^*\}_{i=1}^N$ is the feedback NE for the synthesized game $(A,\{B_i\}_{i=1}^N,\{Q_i^*\}_{i=1}^N,\{R_{ij}\}_{i,j=1}^N)$.\qedsymbol

The next result is the consequence of the previous theoretical results and is important for practical implementation of the algorithm since the results before are valid for infinitely many iterations. 
\begin{thm}
For each iteration $p=0,1,\dots$, there exists a set of learning rates $\{\alpha_i\}_{i=1}^N$ such that $\{K_i^{(p+1)}\}_{i=1}^N$ is the stabilizing solution for \eqref{IOC} and, as a result, the dynamics $A - \sum_{j=1}^N B_j F_j^{(p+1)}$ is stable. 
\end{thm}
\noindent\textbf{Proof}. One can check that $K_i^{(p)}$ linearly affects $F_i^{(p)}$. The initial $K_i^{(p)}$ for $p=0$ is stabilizing as well as the terminal one $K_i^*$ because of \eqref{alggoal}. Hence, referring to \cite{Bertsekas/99}, we know that by choosing an appropriate set of $\{\alpha_i\}_{i=1}^N$ one can always have the next iteration of $K_i^{(p)}$, i.e., $K_i^{(p+1)}$ being a stabilizing solution of \eqref{IOC}. Thus, at each iteration, a game described by $(A,\{B_i\}_{i=1}^N,\{Q_i^{(p+1)}\}_{i=1}^N,\{R_{ij}\}_{i,j=1}^N)$ has an NE feedback $F_i^{(p+1)} = R_{ii} B_i^\top K_i^{(p+1)}$.\qedsymbol

\subsection{Characterization of the Solutions}\label{solcharact}
This section provides a result that allows to adjust the output of \textbf{Algorithm 1}.

Note that we are looking for $\{Q_i^*,R_{ij}\}_{i,j=1}^N$ such that with the dynamics $(A,B_1,\dots,B_N)$ \eqref{AREs} has a stabilizing solution $\{K_i^*\}_{i=1}^N$ satisfying $R_{ii,d}^{-1} B_i^\top K_{i,d} = R_{ii}^{-1} B_i^\top K_i^*$ for $i=1,\dots,N$. Since $R_{ii}>0$, $B_i^\top K_i^* = R_{ii} R^{-1}_{ii,d} B_i^\top K_{i,d}$ for $i=1,\dots,N$. If any of $B_i$ has no full rank, there might be an infinite number of possible $K_i^*$ \cite{lian_robust_2021}.  

\begin{rem}
    All possible outputs of \textbf{Algorithm 1}, i.e.,  $Q_i^*, R_{ij}, K_i^*$, $i,j=1,\dots,N$, satisfy the following equality 
\begin{align}\begin{split}
    & A^\top (K_{i,d} - K_i^*) + (K_{i,d} - K_i^*) A + (Q_{i,d} - Q_i^*) + \\
    & \sum_{j=1}^N F_j^{*\top} (R_{ij,d} - R_{ij}) F_j^* - \\
    & \big(\sum_{j=1}^N F_j^{*\top} B_j^\top \big) (K_{i,d} - K_i^*) - (K_{i,d} - K_i^*) \big(\sum_{j=1}^N B_j F_j^*\big) = 0,\\
    & \text{where}\quad F_i^* = R_{ii}^{-1} B_i^\top K_i^* = F_{i,d},\quad i=1,\dots,N.
\end{split}
\end{align}
\end{rem}
These equations are obtained by the subtraction of \eqref{AREsConv} from \eqref{invpAREs}.
Let us define
\begin{equation}\label{terms}
    \Delta Q_i = Q_i^* - Q_i^\prime, \quad \Delta K_i = K_i^* - K_i^\prime,\quad \Delta R_{ij} = R_{ij} - R_{ij}^\prime,
\end{equation}
where $Q_i^*,R_{ij}$ and $K_i^*$ are the output of \textbf{Algorithm 1} for $i,j = 1,\dots, N$.

\begin{prop}\label{chprop}
Set $\Delta K_i = 0$ and $\Delta R_{ii} = 0$ for $i = 1, \dots, N$ (i.e., $R_{ii} = R_{ii}^\prime$ and $K_i^* = K_i^\prime$). Then, every $Q_i^\prime$ and $R_{ij}^\prime$ for $i,j=1,\dots,N$, $j\neq i$ satisfying 
\begin{equation}
    (Q_i^* - Q_i^\prime) + \sum_{j=1,j\neq i}^N F_j^{*\top} (R_{ij} - R_{ij}^\prime) F_j^* = 0
\end{equation}
together with $R_{ii}^\prime$ and the dynamics $(A,B_1,\dots,B_N)$ form a new game equivalent to $(A,\{B_i\}_{i=1}^N,\{Q_{i,d}\}_{i=1}^N,\{R_{ij,d}\}_{i,j=1}^N)$.
\end{prop}
This is a consequence of a re-scaling of the parameters that does not affect the feedback laws 
\begin{equation}
    F_i^* = F_{i,d} = F_i^\prime = R_{ii}^\prime B_i^\top K_i^\prime.
\end{equation}
Hence, we can adjust $Q_i^\prime$ as 
\begin{equation}\label{qadj}
    Q_i^\prime = Q_i^* + \sum_{j=1,j\neq i}^N F_j^{*\top} (R_{ij} - R_{ij}^\prime) F_j^*
\end{equation}
or $R_{ij}$ for $i\neq j$ in a desired way scaling $Q_i^\prime$. Thus, we can generate an infinite number of possible equivalent games and relax the assumption on definiteness of $R_{ij}$, $i,j=1,\dots,N$, $j\neq i$. 

\section{Model-free Inverse Reinforcement Learning Algorithm}\label{MFsec}
This section present the model-free extension of \textbf{Algorithm 1}. Real-world applications rarely assume the knowledge of the model of the systems. There are three steps in the algorithm presented before that use the system's dynamics matrices - computation of the solution for the initialized game, gradient descent update and the evaluation of the cost function's parameter upgrade. Although there were a number of works dedicated to partially model-free or model-free methods to solve AREs (e.g. \cite{vrabie_adaptive_2009},\cite{vrabie_adaptive_2011}, \cite{vamvoudakis_non-zero_2015},\cite{vamvoudakis_q-learning_2017}), to extend our algorithm, we use the ideas presented in \cite{jiang_computational_2012}, \cite{modares__2015}.

\subsection{Model-free Computation of the Initialized Solution}\label{VA}
After the initialization of the game parameters $R_{ij}$ and $Q_i^{(0)}$, we need to solve the synthesized game. Using the demonstrated trajectories $(x_d,\{u_{i,d}\}_{i=1}^N)$, we use the auxiliary controls
\begin{equation}\label{auxcont}
    u_i^{(k,0)} = - F_i^{(k,0)} x_d,\quad i=1,\dots,N 
\end{equation}
where $k=0,1,\dots$ is the iteration for the step 3 of the algorithm. Using these controls we rewrite the dynamics 
\begin{equation}\label{Acl}
    \dot{x}_d = A x_d + \sum_{i=1}^N B_i u_{i,d} =  A x_d + \sum_{i=1}^N B_i u_i^{(k,0)} + \sum_{i=1}^N B_i (u_{i,d} - u_i^{(k,0)}).
\end{equation}
Using \eqref{auxcont}, we extend the dynamics as
\begin{equation}
    \dot{x}_d = A_{cl}^{(k,0)} x_d + \sum_{i=1}^N B_i (u_{i,d} - u_i^{(k,0)})
\end{equation}
where $A_{cl}^{(k,0)} = A - \sum_{i=1}^N B_i F_i^{(k,0)}$.

Next, for each $i=1,\dots,N$ we multiply \eqref{lyapiter} by $x^\top$ and $x$ to get 
\begin{align}
\begin{split}
    &x_d^\top(A - \sum_{i=1}^N B_i F_i^{(k,0)})^\top K_i^{(k+1,0)}x_d + x_d^\top K_i^{(k+1,0)} (A - \sum_{i=1}^N B_i F_i^{(k,0)}) x_d =\\
    & -x_d^\top (Q_i^{(0)} + \sum_{j=1}^N F_j^{(k,0)\top} R_{ij} F_j^{(k,0)})x_d.
\end{split}
\end{align}
Rewriting the dynamics term, the following equations hold
\begin{align}
    \begin{split}
            & x_d^\top(A_{cl}^{(k,0)} - \sum_{i=1}^N B_i (F_{i,d} - F_i^{(k,0)}))^\top K_i^{(k+1,0)}x_d +\\
            & x_d^\top K_i^{(k+1,0)} (A_{cl}^{(k,0)} - \sum_{i=1}^N B_i (F_{i,d} - F_i^{(k,0)})) x_d =\\
            & -x_d^\top (Q_i^{(0)} + \sum_{j=1}^N F_j^{(k,0)\top} R_{ij} F_j^{(k,0)} + 2 \sum_{j=1}^N (F_{j,d} - F_j^{(k,0)})^\top B_j^\top K_i^{(k+1,0)})x_d.
    \end{split}
\end{align}
Using \eqref{demdyn}, \eqref{auxcont} and \eqref{Acl}, we get
\begin{align}
    \begin{split}
            & \dot{x}_d^\top K_i^{(k+1,0)}x_d + x_d^\top K_i^{(k+1,0)}\dot{x}_d =x_d^\top ( - Q_i^{(0)}  - \sum_{j=1}^N F_j^{(k,0)\top} R_{ij} F_j^{(k,0)} - \\ 
            & 2 (F_{i,d} + F_i^{(k,0)})^\top B_i  K_i^{(k+1,0)} -
            2  \sum_{j\neq i}^N (F_{j,d} - F_j^{(k,0)})^\top Y_{ji}^{(k+1)} ) x_d.
    \end{split}
\end{align}
where $Y_{ji}^{(k+1)} = B_j^\top K_i^{(k+1,0)}$ is another auxiliary variable for $j\neq i, i=1,\dots,N$. Considering \eqref{trgtfdbkl} and following the ideas of \cite{jiang_computational_2012} and \cite{modares__2015}, we integrate the above equation from $t$ to $t+T$ as follows 
\begin{align}\label{mfinitialsol}
    \begin{split}
        & x_d^\top (t+T) K_i^{(k+1,0)} x_d (t+T) - x_d^\top (t) K_i^{(k+1,0)} x_d (t) - \\ 
        & 2 \int_{t}^{t+T} (u_{i,d} + F_i^{(k,0)}x_d)^\top R_{ii}  F_i^{(k+1,0)} x_d \,d\tau - \\
        & 2 \sum_{j\neq i}^N \int_{t}^{t+T} (u_{j,d} + F_j^{(k,0)}x_d))^\top Y_{ji}^{(k+1)} x_d\,d\tau = \\
        & - \int_{t}^{t+T} x_d^\top( Q_i^{(0)} + \sum_{j=1}^NF_j^{(k,0)\top} R_{ij} F_j^{(k,0)}) x_d\,d\tau,\quad i=1,\dots,N.
    \end{split}
\end{align}

Let us consider the above for $k=0$. Set the initial stabilizable feedback laws $F_i^{(0,0)} = \hat{F}_i$ where $\hat{F}_i$ is estimated in step 2 for $i=1,\dots,N$. Then, the unknowns in the above equations are 
\begin{equation}
    K_i^{(k+1,0)},\quad F_i^{(k+1,0)},\quad Y_{ji}^{(k+1)},\quad i,j=1,\dots,N,\, i\neq j
\end{equation}
and each of them is built using $K_i^{(k+1)}$ and $K_j^{(k+1)}$ for $i$ and $j\neq i$. Thus, we solve \eqref{mfinitialsol} with respect to the mentioned unknowns till $\lVert K_i^{(k+1,0)} - K_i^{(k,0)}\rVert \leq \epsilon_i$ where $\epsilon_i > 0$ is a small constant that describes a measure of precision for $i=1,\dots,N$.

To perform the next steps, matrices $B_i$, $i=1,\dots, N$, are needed. One way to evaluate these matrices is to use the computed values of $K_i^{(k+1)}$, $F_i^{(k+1),0}$ and $Y_{ji}^{(k+1)}$. Recall that 
\begin{align}
    \begin{split}
        F_i^{(k+1,0)} &= R_{ii}^{-1} B_i^\top K_i^{(k+1,0)},\\
        Y_{ji}^{(k+1)} &= B_j^\top K_i^{(k+1,0)}.
    \end{split}
\end{align}
Using the computed values from equation \eqref{mfinitialsol} associated with any player $i\in\{1,\dots,N\}$, the control input matrices can be evaluated as 
\begin{align}
    \begin{split}
        B_i &= (R_{ii} F_i^{(k+1,0)} (K_i^{(k+1,0)})^{-1})^\top,\\
        B_j &= (Y_{ji}^{(k+1)} (K_i^{(k+1,0)})^{-1})^\top,\quad j\neq i.
    \end{split}
\end{align}
Note that the inverse $K_i^{(k+1,0)}$ exists because the initialized cost function parameters guarantee $K_i^{(k+1,0)} > 0$ for $i=1,\dots,N$.

\begin{thm}\label{thmequiv}
    The solution $\{K_i^{(k+1,0)}\}_{i=1}^N$ of \eqref{mfinitialsol} is a unique positive definite stabilizing solution and is the same as the solution of \eqref{mbinitialsol}. 
\end{thm}
\noindent\textbf{Proof}. We give a short proof here that follows \cite{jiang_computational_2012} and \cite{modares__2015}. We can reverse engineer \eqref{mfinitialsol} taking its $\lim_{T\to0}$ and using L'Hopital's rule \cite{rudin1976principles} to derive \eqref{mbinitialsol}. According to \cite{li_1995}, the solution of \eqref{mbinitialsol} is a unique positive definite solution for the cost function parameters satisfying $Q_i>0,R_{ii}>0,R_{ij}\geq0$ for $i\neq j$, $i,j=1,\dots,N$. Thus, we conclude that \eqref{mfinitialsol} has the same solution as \eqref{mbinitialsol} that is a stabilizing positive definite one.\qedsymbol

After the computation of the initial solution $\{K_i^{(k+1,0)}\}_{i=1}^N$  is accomplished , as it is done in step 4, we drop the iteration counter and set 
\begin{equation}\label{initmf}
    K_i^{(k+1,0)} = K_i^{(0)},\quad F_i^{(k+1,0)} = F_i^{(0)},\quad i=1,\dots,N.
\end{equation}

\subsection{Mode-free Inverse Update of the Parameters}
Since we evaluated $B_i$ Step $5$ can be used as it is in Algorithm $1$.
\begin{rem}
    In fact, using \eqref{eqfdblw} and the values in \eqref{initmf}, one can conclude the following
\begin{align}
    \begin{split}
        & F_i^{(0)} = R_{ii}^{-1} B_i^\top K_i^{(0)},\\
        & F_i^{(0)} (K_i^{(0)})^{-1} = R_{ii}^{-1} B_i^\top
    \end{split}
\end{align}
because $K_i^{(0)}$ is guaranteed to be a positive definite solution of \eqref{mfinitialsol} as it is shown in Theorem \ref{thmequiv} for $i=1,\dots,N$. Thus, step 5 can also be rewritten as  
    \begin{align}\label{mfupdate}
    \begin{split}
	        &K_i^{(p+1)} = K_i^{(p)} - \alpha_i \Big(d^{(p)\top}_i  F_i^{(0)} (K_i^{(0)})^{-1} + (K_i^{(0)})^{-1} F_i^{(0)\top} d_i^{(p)}\Big),\\
            & F^{(p+1)}_i = F_i^{(0)} (K_i^{(0)})^{-1} K_i^{(p+1)}. 
	 \end{split}   
     \end{align}
\end{rem} 

The last update, step $5$ in \eqref{IOC}, can also modified to avoid using the unknown matrices. Following the approach used in \ref{VA}, one can rewrite \eqref{IOC} as 
\begin{align}
     \begin{split}
        & x_d^\top Q_i^{(p+1)} x_d =  x_d^\top (-\sum_{j=1}^N F_j^{(p+1)\top} R_{ij} F_j^{(p+1)} - \\
        & (A_{cl}^{(p+1)} - \sum_{j=1}^N B_j (F_{j,d} - F_j^{(p+1)}))^\top K_i^{(p+1)} - \\
        &  K_i^{(p+1)} (A_{cl}^{(p+1)} - \sum_{j=1}^N B_j (F_{j,d} - F_j^{(p+1)})) + \\
        & 2 \sum_{j=1}^N (F_{j,d} - F_j^{(p+1)}))^\top B_j^\top K_i^{(p+1)}) x_d .
    \end{split}
\end{align}
Integrating both sides of the above equation from $t$ to $t+T^\prime$, we get 
\begin{align}\label{IOCmf}
    \begin{split}
        & \int_{t}^{t+T^\prime} x_d^\top Q_i^{(p+1)} x_d\,d\tau = - \int_{t}^{t+T^\prime} x_d^\top \sum_{j=1}^N F_j^{(p+1)\top} R_{ij} F_j^{(p+1)} x_d\,d\tau -\\
        & x_d^\top(t+T^\prime) K_i^{(p+1)} x_d(t+T^\prime) + x_d^\top(t) K_i^{(p+1)} x_d(t) - \\
        & 2 \int_{t}^{t+T^\prime}  \sum_{j=1}^N  (u_{j,d} + F_j^{(p+1)}x_d)^\top B_j^\top K_i^{(p+1)} x_d.
    \end{split}
\end{align}
Since \eqref{mfupdate} provides us $K_i^{(p+1)}$,$F_i^{(p+1)}$ and the trajectories $(x_d,\{u_{i,d}\}_{i=1}^N)$ are given, $Q_i^{(p+1)}$ can be evaluated. The way it can be done is shown in the next section. All the steps for the model-free \textbf{Algorithm 2} are shown below.

\begin{algorithm}\label{alg2}
	\caption{Model-free Inverse Reinforcement Learning Algorithm}
	\begin{enumerate}
	    \item Initialize $R_{ii}> 0$, $R_{ij}\geq 0$ and $Q_i^{(0)}> 0$ for $i,j = 1,\dots,N$, $i\neq j$. Sample data from demonstrated $(x,\{u_{i,d}\}_{i=1}^N)$ to generate $(\hat{x},\{\hat{u}_{i,d}\}_{i=1}^N)$. Set $k=0$ and $p=0$.
	    \item Derive estimation of  $F_{i,d}$ using the sampled data as
	    \begin{equation}
	        \hat{F}_i =  -\hat{u}_{i,d}\hat{x}_d^\top(\hat{x}_d \hat{x}_d^\top)^{-1}.
	    \end{equation}
	    \item Set $F_i^{(0,0)} = \hat{F}_i$ for $i=1,\dots,N$, solve \eqref{mfinitialsol} with respect to $K_i^{(k+1)}$, $F_i^{(k+1,0)}$ and $Y_{ji}^{(k+1)}$ for $i,j=1,\dots,N$, $j\neq i$. Compute $B_i$ for $i=1,\dots,N$. Set $k = k+1$ till $\lVert K_i^{(k+1,0)} - K_i^{(k,0)} \rVert < \epsilon_i$ where $\epsilon_i$ is a small positive constant for $i=1,\dots, N$. 
 	    \item Set $K_i^{(0)} = K_i^{(k+1,0)}$, $F_i^{(0)} = F_i^{(k+1,0)}$. Compute $F_i^{(0)} (K_i^{(0)})^{-1}$
	    and evaluate the difference
	    \begin{equation}
	    d^{(p)}_i = F_i^{(p)} - \hat{F}_i.
	    \end{equation}
	    \item Update $K^{(p+1)}_i$ and $F_i^{(p+1)}$ for $i=1,\dots,N$ as in \eqref{mfupdate}.
	    \item Perform evaluation of $Q_i^{(p+1)}$ from \eqref{IOCmf}.       
	    \item Set $p = p + 1$. Perform steps $4$-$6$ till $\tr(d_i^{(p)\top} d_i^{(p)}) < \delta_i$ where $\delta_i$ is a small positive constant for $i=1,\dots,N$.
	\end{enumerate}
\end{algorithm}
\begin{rem}\label{step6}
As for \textbf{Algorithm 1}, the implementation of \textbf{Algorithm 2} does not necessarily require the iterative update of $Q_i^{(p+1)}$ in step 6. This update might be done only \emph{once} after the desired precision $\delta_i$ is achieved, i.e., after getting $K_i^{(p+1)}$ in step 5 such that $\tr(d_i^{(p)\top} d_i^{(p)}) < \delta_i$ for $i=1,\dots,N$. 
\end{rem}
\subsection{Implementation of the algorithm}
In this section, we show \emph{one} possible way to implement \textbf{Algorithm 2} which is partially based on \cite{jiang_computational_2012}. For other ways to use the proposed algorithm, the reader can check \cite{modares__2015} ,\cite{vamvoudakis_non-zero_2015},\cite{xue_inverse_2021}. To avoid any confusion due to indexes and terms, we show the algorithm implementation for the two-player case, i.e., $N=2$. We hope the below description of the implementation clarifies for the reader the implementation of the algorithm in the multiplayer case. 

Firstly, we show how to perform evaluation of $K_i^{(k+1,0)}$, $F_i^{(k+1,0)}$ and $Y_{ji}^{(k+1,0)}$ in step $3$ from \eqref{mfinitialsol}. Following \cite{jiang_computational_2012}, the following notations are introduced
\begin{align}
    \begin{split}
        \hat{K}_i &= [k_{i,11},2k_{i,12},\dots,2k_{i,1n}, k_{i,22},2k_{i,23},\dots,k_{i,nn}]^\top\in\mathbb{R}^{n(n+1)/2},\\
        \hat{x} &= [x_1^2,x_1 x_2, \dots, x_1 x_n, x_2^2, x_2 x_3,\dots, x_n^2]^\top\in\mathbb{R}^{n(n+1)/2}.
    \end{split}
\end{align}
where $k_{i,l_1 l_2}$ is a particular element of matrix $K_i$, i.e., $(K_i)_{l_1l_2}$ for $l_1,l_2=1,\dots,n$. We use the following property of the Kronecker product 
\begin{equation}
    (c^\top \otimes a^\top) \text{vec}(B) = a^\top B c.
\end{equation}
Thus, one can rewrite terms in \eqref{mfinitialsol} as 
\begin{align}
    \begin{split}
        & x_d^\top (t+T) K_i^{(k+1,0)} x_d (t+T) - x_d^\top (t) K_i^{(k+1,0)} x(t) =\\
        & (\hat{x}(t+T) - \hat{x}(t)) \hat{K}_i^{(k+1,0)}, \\
        &  (u_{i,d} + F_i^{(k,0)}x_d)^\top R_{ii} F_i^{(k+1,0)} x_d= ((x_d^\top\otimes u_{i,d}^\top) (I_n\times R_{ii}) + \\& (x_d\otimes x_d)(I_n\otimes F_i^{(k,0)\top}R_{ii}))\text{vec}(F_i^{(k+1,0)}),\\
        & (u_{j,d} + F_j^{(k,0)}x_d)) Y_{ji}^{(k+1)} x_d = ((x_d^\top\otimes u_{j,d}^\top)  + \\& (x_d\otimes x_d)(I_n\otimes F_j^{(k,0)\top}))\text{vec}(Y_{ji}^{(k+1)}).\\
    \end{split}
\end{align}
In addition to the above, we define $\delta_{xx}$, $I_{xx}$ and $I_{xu_i}$
as 
\begin{align}\label{nottttations}
    \begin{split}
        \delta_{xx}&= [\hat{x}(t_1)-\hat{x}(t_0),\hat{x}(t_2) - \hat{x}(t_1),\dots, \hat{x}(t_s) - \hat{x}(t_{s-1})]^\top,\\
        I_{xx} &= \int_{t_0}^{t_1} (x_d\otimes x_d)\,d\tau,\int_{t_1}^{t_2} (x_d\otimes x_d)\,d\tau,\dots,\int_{t_{s-1}}^{t_s} (x_d\otimes x_d)\,d\tau]^\top \\
        I_{xu_i} &= \int_{t_0}^{t_1} (x_d\otimes u_{i,d})\,d\tau,\int_{t_1}^{t_2} (x_d\otimes u_{i,d})\,d\tau,\dots,\int_{t_{s-1)}}^{t_s} (x_d\otimes u_{i,d})\,d\tau]^\top\\ 
    \end{split}
\end{align}
where $0\leq t_{l-1}\leq t_{l}$ for $l\in\{0,1,\dots,s\}$. Although the data intervals do not need to be equal, in our simulation presented further, we use $t_l-t_{l-1} = T$ for $l\in\{0,1,\dots,s\}$.

Then, \eqref{mfinitialsol} can be rewritten as 
\begin{equation}\label{lq11}
    H_i^{(k)} \begin{pmatrix}
        \hat{K}_i^{(k+1,0)} \\ \text{vec}(F_i^{(k+1,0)}) \\ \text{vec}(Y_{ji}^{(k+1)})\\
    \end{pmatrix}= \Xi_i^{(k)}
\end{equation}
where 
\begin{align}
    \begin{split}
        H_i^{(k)} = &[\delta_{xx}, -2I_{xu_i}(I_n\otimes R_{ii}) - I_{xx}(I_n\otimes F_i^{(k,0)^\top} R_{ii}, \\ 
        & -2I_{xu_i} - I_{xx}(I_n\otimes F_i^{(k,0)^\top}],\\
        \Xi_i^{(k)} = & -I_{xx} (Q_i^{(0)} + \sum_{j=1}^2 F_j^{(k,0)\top} R_{ij} F_j^{(k,0)}).
    \end{split}
\end{align}
Then, \eqref{lq11} can be solved as 
\begin{equation}\label{lq111}
    \begin{pmatrix}
        \hat{K}_i^{(k+1,0)} \\ \text{vec}(F_i^{(k+1,0)}) \\ \text{vec}(Y_{ji}^{(k+1)})\\
    \end{pmatrix} = (H_i^{(k)\top} H_i^{(k)})^{-1} H_i^{(k)\top} \Xi_i^{(k)}.
\end{equation}
The equation is solved until the convergence of $\hat{K}_i^{(k+1,0)}$ from which one can recover $K_i^{(k+1,0)}$. Note that the vector of unknowns has $n(n+1)/2 + m_i n + m_j n$ parameters. Thus, we need enough data to satisfy $s \geq n(n+1)/2 + m_i n + m_j n$. 
\begin{rem}
    If $H_i^{(k)}$ is an invertible square matrix, right side of \eqref{lq111} can be computed as ($H_i^{(k)})^{-1} \Xi_i^{(k)}$.
\end{rem}

Although step $6$ of \textbf{Algorithm 2} can be implemented interactively for every new $K_i^{(p+1)}$, one can implement it only once, as suggested in Remark \ref{step6} after the feedback laws converged as a result of the gradient updates, i.e., $\tr(d_i^{(p)\top} d_i^{(p)}) < \delta_i$ for $i=1,\dots,N$. Then, set $K_i^{(p+1)} = K_i^*$ and $F_i^{p+1} = F_i^*$. For that one use the same data as in \eqref{lq11}. We rewrite one of the terms in \eqref{IOCmf} as 
\begin{align}
    \begin{split}
        & (u_{j,d} + F_j^* x_d)^\top B_j^\top K_i^* x_d = \\
        & ((x_d^\top\otimes u_{j,d}) + (x_d\otimes x_d)) \text{vec}(F_j^{*\top} B_j^\top K_i^*),\\
        & x_d^\top Q_i^* x_d = I_{xx}\text{vec}(Q_i^*).
    \end{split}
\end{align}
In addition to the above, we define $\hat{Q}_i^*$ and $I_{qx}$ as
\begin{align}
    \begin{split}
        \hat{Q}_i^* &= [q_{i,11},2q_{i,12},\dots,2q_{i,1n}, q_{i,22},2q_{i,23},\dots,q_{i,nn}]^\top,\\
        I_{qx} &= [\int_{t_0}^{t_1} \hat{x}\,d\tau,\int_{t_1}^{t_2} \hat{x}\,d\tau,\dots,\int_{t_{s-1}}^{t_s} \hat{x}\,d\tau]^\top .
    \end{split}
\end{align}

Then, using \eqref{nottttations}, \eqref{IOCmf} can be rewritten as
\begin{align}\label{lq22}
    I_{qx} \hat{Q}_i^* = \Omega_i 
\end{align}
where 
\begin{align}
\begin{split}
    \Omega_i =&   -I_{xx} \sum_{j=1}^2 \text{vec}(F_j^{*\top} R_{ij} F_j^*) - \\
    & \delta_{xx}\hat{K_i}^* - 2\sum_{j=1}^N (I_{xu_j} + I_{xx}) \text{vec}(F_j^{*\top} B_j^\top K_i^*). 
    \end{split}
\end{align}
Then, \eqref{lq22} can be solved as 
\begin{equation}
    \hat{Q}_i^*= (I_{qx}^\top I_{qx})^{-1} I_{qx}^\top \Omega_i.
\end{equation}
Note, \eqref{lq22} has less unknown parameters than \eqref{lq11} because $\hat{Q}_i^*\in\mathbb{R}^{n(n+1)/2}$. Thus, the previous restriction on $s$ is enough, i.e., $s\geq n(n+1)/2 + m_i n + m_j n$.


\begin{rem}
Proposition \ref{chprop} also valid in the model-free case. Thus, the value of the output of \textbf{Algorithm 2} $Q_i^*$ can be adjusted or the restrictions on definiteness of $R_{ij}$ can be relaxed for $i\neq j$, $i,j =1,\dots,N$.    
\end{rem}

\begin{rem}
Since the equations \eqref{lq11} and \eqref{lq22} are solved as LQ problems, the probing noise should be injected to satisfy persistence of excitation (PE) condition \cite{modares__2015}, \cite{jiang_computational_2012}, \cite{vrabie_adaptive_2011}, \cite{xue_inverse_2021}. The noise can be sinusoids of different frequencies or some random noise. We refer the reader to \cite{ioannou2006} for more details on that matter.  
\end{rem}
Thus, we need to make the following assumption
\begin{assump}\label{noise}
One of the following is true
\begin{itemize}
    \item One can use the estimated stabilizable feedback law $\hat{F}_i$ from \eqref{LQest} to apply control inputs $\hat{u}_i = -\hat{F}_i x + \omega_i(t)$, where $\omega_i(t)$ is a noise term, for $i=1,\dots,N$ to the system for data collection on the range $(t,t_{\bar{N}})$ at $\bar{N}\geq \max(\bar{n},\bar{m})$ points. The collection of additional data is performed once.    

    \item The demonstrated trajectories were generated under the control inputs $u_{i,d} = -F_{i,d} x_d + \omega_i(t)$ where $\omega_i(t)$ is an exponentially decaying noise such that \eqref{lq11} and \eqref{lq22} have a solution. In other words, when the noise decayed significantly, the demonstrated trajectory is $u_{i,d} \approx - F_{i,d} x_d$ for $i=1,\dots,N$. 
    \end{itemize}
\end{assump}
\section{Simulations}\label{SIMsec}
In this section, we present the simulation results of the algorithms developed in this paper. 

\subsection{Model-based Algorithm Simulation}
\begin{figure} 
    \centering
  \subfloat[\label{1a}]{%
       \includegraphics[width=1\linewidth]{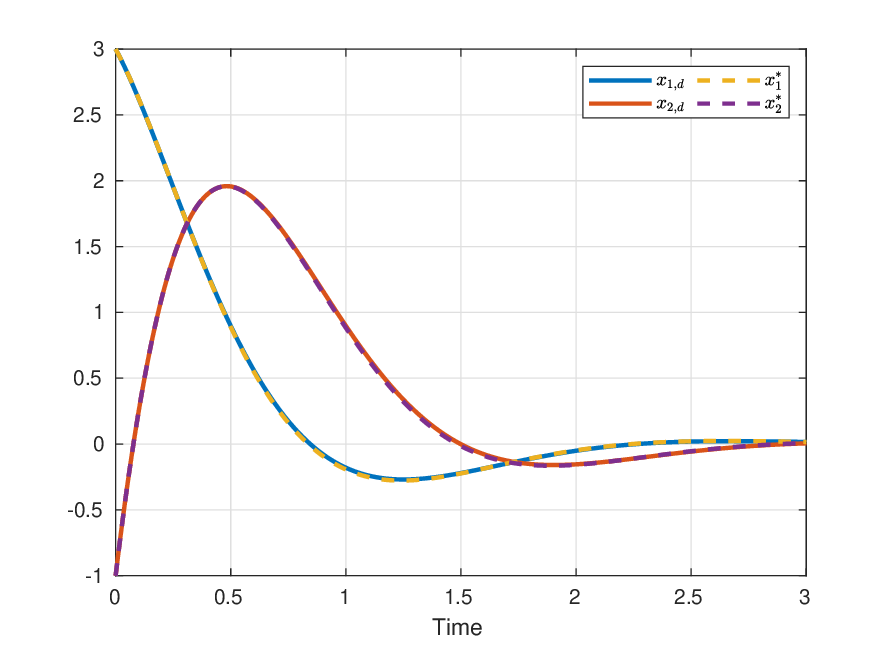}}\\
  \subfloat[\label{1b}]{%
        \includegraphics[width=0.49\linewidth]{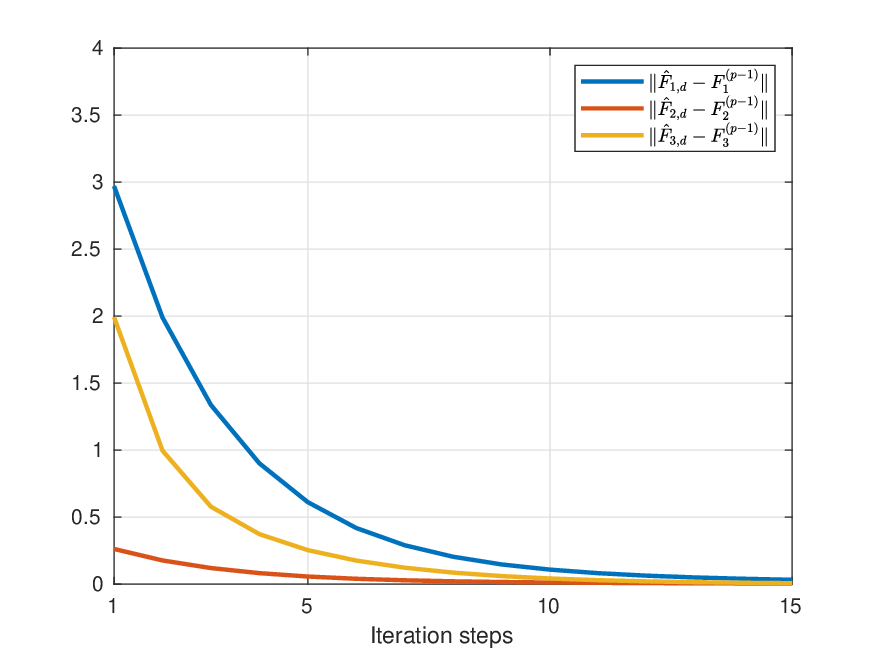}}
    \hfill
  \subfloat[\label{1c}]{%
        \includegraphics[width=0.49\linewidth]{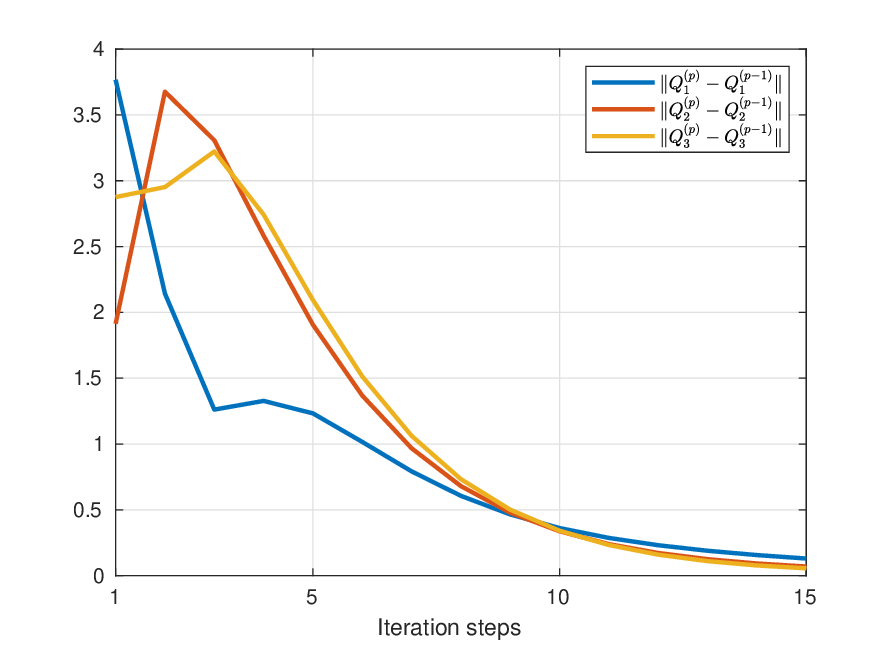}}
  \caption{Algorithm 1: (a) the stability of the demonstrated and resulting dynamics; (b,c) convergence of the norm for iterations of $F_i^{(p)}$ and $Q_i^{(p)}$, respectively.}
  \label{fig1} 
\end{figure}

Consider the following continuous time system dynamics
\begin{equation}
	\dot{x} = A x + \sum_{i=1}^3 B_i u_i, 
\end{equation}
where 
\begin{equation}
	A = \begin{pmatrix}
		3 & -2\\ 4 & -1 	\end{pmatrix},\quad B_1 = \begin{pmatrix}
		1\\ 0
	\end{pmatrix}, \quad B_2 = \begin{pmatrix}0 \\ 1\end{pmatrix},\quad B_3 = \begin{pmatrix}
	    1\\1
	\end{pmatrix}.
\end{equation}
The demonstrated NE trajectories are generated for the game with the following weight matrices
\begin{align}
\begin{split}	
 &Q_{1,d} = \begin{pmatrix}
		7 & 2\\ 2 & 5
	\end{pmatrix},\quad Q_{2,d} = 3 I_{2\times 2}, \quad Q_{3,d} = I_{2\times 2},\\
    &R_{11,d} = 3,\quad R_{12,d} = 1, \quad R_{13,d} = 1,\\
    &R_{21,d} = 1,\quad R_{22,d} = 2, \quad R_{23,d} = 0,\\
    &R_{31,d} = 0,\quad R_{32,d} = 1, \quad R_{33,d} = 4.
\end{split}
\end{align}
Given this game, $F_{1,d}$, $F_{2,d}$ and $F_{3,d}$ are 
\begin{align}
    \begin{split}
        & F_{1,d} = \begin{pmatrix} 4.2499 & -0.9409 \end{pmatrix}, \\
        & F_{2,d} = \begin{pmatrix} -0.4108 & 0.9187 \end{pmatrix}, \\
        & F_{3,d} = \begin{pmatrix} 0.2334 & 0.1295 \end{pmatrix}, 
    \end{split}
\end{align}
with the symmetric solution of AREs
\begin{align}
\begin{split}
    K_{1,d} &= \begin{pmatrix}
        12.7497 & -2.8228 \\ -2.8228 & 3.7172
    \end{pmatrix},\\ 
    K_{2,d} &= \begin{pmatrix}
        4.8994 & -0.8216 \\ 0.8216 & 1.8373
    \end{pmatrix},\\ 
    K_{3,d} &= \begin{pmatrix}
        0.8116 & 0.1222 \\ 0.1222 & 0.3956
    \end{pmatrix}.
\end{split}
\end{align}
The initialized parameters are the following
\begin{align}
\begin{split}	
 &Q_1^{(0)} = I_{2\times 2},\quad Q_2^{(0)} = I_{2\times 2}, \quad Q_3^{(0)} = I_{2\times 2},\\
    &R_{11} = 3,\quad R_{12} = 2, \quad R_{13} = 1,\\
    &R_{21} = 2,\quad R_{22} = 3, \quad R_{23} = 1,\\
    &R_{31} = 2,\quad R_{32} = 3, \quad R_{33} = 1.
\end{split}
\end{align}
The learning rates are set to $\alpha_1 = 1.5$, $\alpha_2 = 1.5$, $\alpha_3=0.15$.\\
The solution generated by the algorithm is 
\begin{align}
\begin{split}
    Q_1^* &= \begin{pmatrix}
        6.2118 & 9.0007 \\ 9.0007 & -1.9440
    \end{pmatrix},\\ 
    Q_2^* &= \begin{pmatrix}
        -15.6574 & -2.6946 \\ -2.6946 & 4.0308
    \end{pmatrix},\\ 
    Q_3^* &= \begin{pmatrix}
        -13.3547 & -4.2811 \\ -4.2811 & -0.3872
    \end{pmatrix}.
\end{split}
\end{align}
with 
\begin{align}
    \begin{split}
        & F_1^* = \begin{pmatrix} 4.2398 & -0.9103 \end{pmatrix}, \\
        & F_2^* = \begin{pmatrix} -0.4149 & 0.9178 \end{pmatrix}, \\
        & F_3^* = \begin{pmatrix} 0.2384 & 0.1245 \end{pmatrix}, 
    \end{split}
\end{align}
and the symmetric solution of AREs given by
\begin{align}
\begin{split}
    K_1^* &= \begin{pmatrix}
        12.7925 & -2.7461 \\ -2.7461 & 2.1820
    \end{pmatrix},\\ 
    K_2^* &= \begin{pmatrix}
        3.5608 & -1.2426 \\ -1.2426 & 2.7543
    \end{pmatrix},\\ 
    K_3^* &= \begin{pmatrix}
        2.2276 & -1.9906 \\ -1.9906 & 2.1166
    \end{pmatrix}.
\end{split}
\end{align}
The resulting dynamics $A + \sum_{i=1}^3 B_i F_i^*$ are stable as shown in Figure \ref{1a}. The convergence of the iterative procedure is shown in Figures \ref{1b} and \ref{1c}.

\begin{rem}
    The reader might notice that the learning rate for players $1,2$ and player $3$ differ. The reason is that for $\alpha_3=\alpha_1=\alpha_2$ the overshooting of the gradient descent method is observed. In fact, an adaptive learning rate might be used, e.g. Polyak step-size and the line search method \cite{sun2006optimization}.
\end{rem}

\subsection{Model-free Algorithm Simulation}

\begin{figure} 
    \centering
  \subfloat[\label{2a}]{%
       \includegraphics[width=1\linewidth]{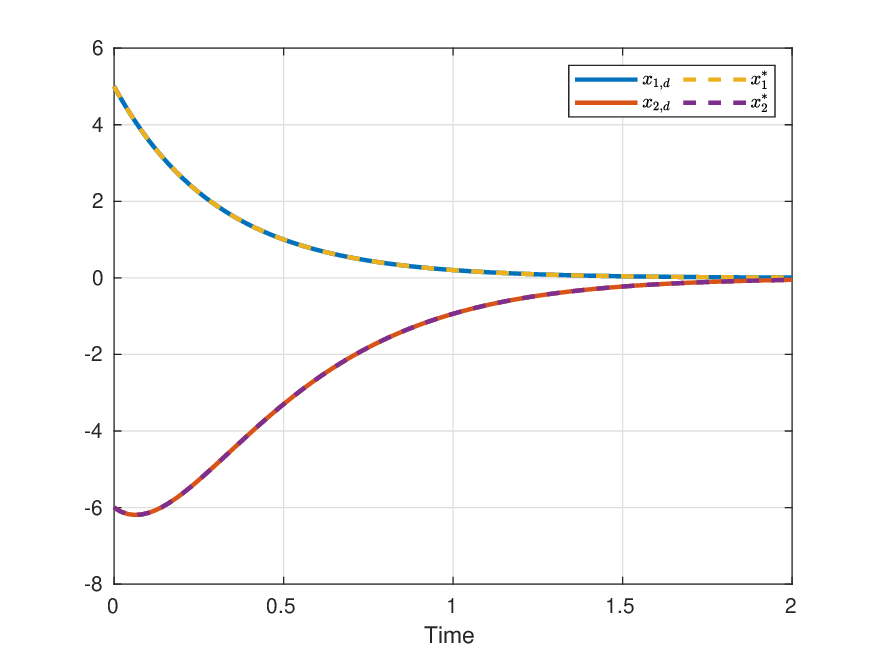}}\\
  \subfloat[\label{2b}]{%
        \includegraphics[width=0.49\linewidth]{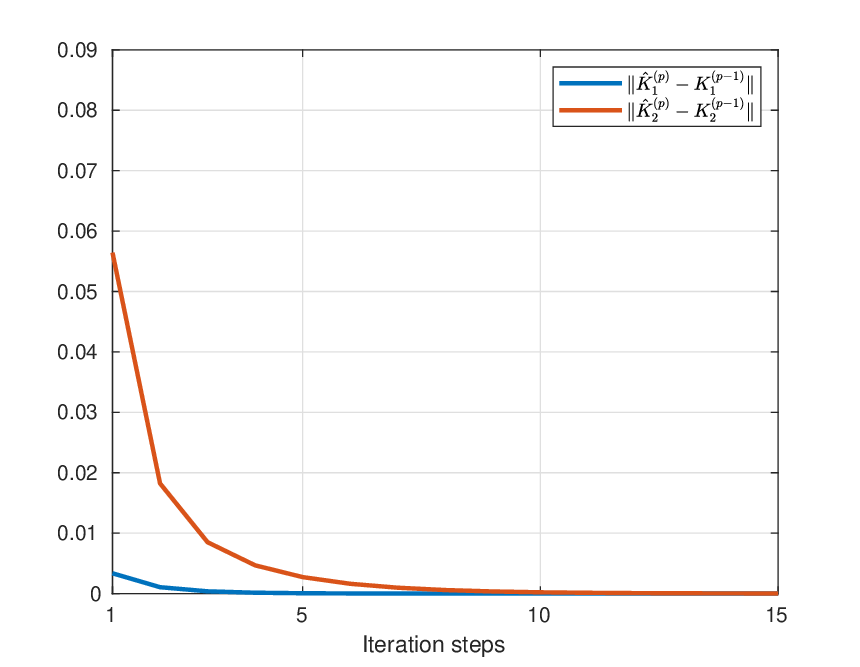}}
    \hfill
  \subfloat[\label{2c}]{%
        \includegraphics[width=0.49\linewidth]{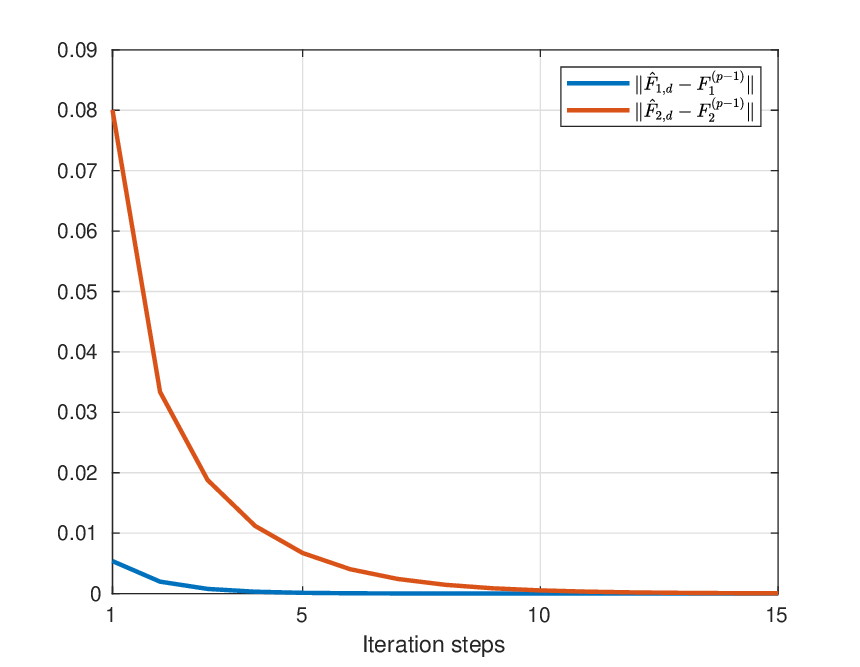}}
  \caption{Algorithm 2: (a) the stability of the observed and resulting dynamics; (b,c) convergence of the norm for iterations of $K_i^{(p)}$ and $F_i^{(p)}$, respectively.}
  \label{fig2} 
\end{figure}

Consider the following continuous time system dynamics
\begin{equation}
	\dot{x} = A x + \sum_{i=1}^2 B_i u_i, 
\end{equation}
where 
\begin{equation}
	A = \begin{pmatrix}
		3 & 0\\ 0 & -4 	\end{pmatrix},\quad B_1 = \begin{pmatrix}
		1\\ 1
	\end{pmatrix}, \quad B_2 = \begin{pmatrix}0 \\ 1\end{pmatrix}.
\end{equation}
The demonstrated NE trajectories are generated for the game with the following weight matrices
\begin{align}
\begin{split}	
 &Q_{1,d} =2 I_2 \quad Q_{2,d} = 3 I_2,\\
    &R_{11,d} = 2,\quad R_{12,d} = 1, \\
    &R_{21,d} = 1,\quad R_{22,d} = 6. \\
\end{split}
\end{align}
Given this game, $F_{1,d}$ and $F_{2,d}$ are 
\begin{align}
    \begin{split}
        & F_{1,d} = \begin{pmatrix} 6.2586 & 0.0186 \end{pmatrix}, \\
        & F_{2,d} = \begin{pmatrix} -0.0532 & 0.0620 \end{pmatrix}, 
    \end{split}
\end{align}
with the symmetric solution of AREs
\begin{align}
\begin{split}
    K_{1,d} &= \begin{pmatrix}
        12.7267 & -0.2095 \\ -0.2095 & 0.2466
    \end{pmatrix},\\ 
    K_{2,d} &= \begin{pmatrix}
        7.0811 & -0.3192 \\ -0.3192 & 0.3719
    \end{pmatrix}.
\end{split}
\end{align}
Firstly, given the demonstrated trajectories of the game described above, we estimate \eqref{LQest} $\hat{F}_1,\hat{F}_2$. Then, following Assumption \ref{noise}, for additional data collection we applied the following controller 
\begin{equation}
\hat{u}_i = - \hat{F}_i x + \omega_i
\end{equation}
for $i=1,2$ where $\omega_i (t) =  100\sum_{k=1}^{100} \sin(c_k t)$ and $c_k$ for $k=1,\dots,100$ is a random number selected in the range $[-500,500]$ \cite{jiang_computational_2012}. Data are collected at $0.01$ sec during $2$ seconds. 
Then, using the collected data and the initialized parameters below
\begin{align}
\begin{split}	
 &Q_1^{(0)} = I_2,\quad Q_2^{(0)} = I_2,\\
    &R_{11} = 3,\quad R_{12} = 0,\\
    &R_{21} = 0,\quad R_{22} = 3,
\end{split}
\end{align}
we derive solution for the initialized game as 
\begin{align}
    \begin{split}
        & K_1^{(k+1,0)} = \begin{pmatrix}
            6.3546 & -0.1011 \\ -0.1011 & 0.1212
        \end{pmatrix},\\
        & K_2^{(k+1,0)} = \begin{pmatrix}
            6.3538 & -0.1050 \\ -0.1050 & 0.1230
        \end{pmatrix},
    \end{split}
\end{align}
with the following equilibrium feedback laws 
\begin{align}
    \begin{split}
        & F_1^{(k+1,0)} = \begin{pmatrix} 6.2535 & 0.0202 \end{pmatrix}, \\
        & F_2^{(k+1,0)} = \begin{pmatrix} -0.1050 & 0.1230 \end{pmatrix}. 
    \end{split}
\end{align}
The learning rates are set to $\alpha_1 = 0.3$,  $\alpha_2 = 0.4$.\\
As suggested in Remark \ref{step6}, we perform  \eqref{lq111} only once after getting convergence of $F_i$ for $i=1,2$. The solution generated by the algorithm is 
\begin{align}
\begin{split}
    Q_1^* &= \begin{pmatrix}
        1.0284 & 0.0034 \\ 0.0034 & 0.9648
    \end{pmatrix},\\ 
    Q_2^* &= \begin{pmatrix}
        1.6420 & 0.0039 \\ 0.0039 & 0.4998
    \end{pmatrix}.
\end{split}
\end{align}
with 
\begin{align}\label{fsol}
    \begin{split}
        & F_1^* = \begin{pmatrix} 6.2586 & 0.0186 \end{pmatrix}, \\
        & F_2^* = \begin{pmatrix} -0.0532 & 0.0620 \end{pmatrix},
    \end{split}
\end{align}
and the symmetric solution of AREs given by
\begin{align}\label{ksol}
\begin{split}
    K_1^* &= \begin{pmatrix}
        6.3588 & -0.1002 \\ -0.1002 & 0.1187
    \end{pmatrix},\\ 
    K_2^* &= \begin{pmatrix}
        0.3537 & -0.0532 \\ -0.0532 & 0.0620
    \end{pmatrix}.
\end{split}
\end{align}
The resulting dynamics $A + \sum_{i=1}^2 B_i F_i^*$ is stable, as shown in Figure \ref{2a}. The convergence of the iterative procedure is shown in Figures \ref{2b} and \ref{2c}.
\begin{rem}
As suggested in the solution characterization section \eqref{qadj}, one can change the algorithm output, preserving the game equivalence. For example, set a new $R^\prime_{21} = -1$ instead of $R_{21} = 0$ used as initialized parameter, relaxing the positive definiteness assumption on $R_{21}$. Then, the game with the same parameters as above, except 
\begin{equation}
    Q^\prime_2 = Q_2^* + F_1^{*\top} (R_{21} - R_{21}^\prime) F_1^* = \begin{pmatrix}
        40.8124 & 0.1200 \\ 0.1200 & 0.5002
    \end{pmatrix}
\end{equation}
instead of $Q_2^*$ and $R_{21}^\prime=-1$ instead of $R_{21}=0$ is also equivalent to the observed game, i.e., it has solution given by \eqref{ksol} and \eqref{fsol}.
\end{rem}

\section{Conclusion}\label{CONsec}
In this paper, we provide algorithms to solve the inverse problem for linear-quadratic nonzero-sum differential games. Both model-based and model-free versions were introduced. We showed that the algorithms' output is the set of weight matrices that together with the dynamics matrices form an equivalent game for one of the players. After showing the convergence of the algorithms to a desired output, we also provided solution characterizations and showed how the algorithms' output could be adjusted. The effectiveness of the algorithm was demonstrated via simulations. We discussed how the algorithms could be implemented with low (as much as possible) computational cost. The presented algorithms can be extended for the case of non-linear dynamics of the form $f(x)+\sum_{i=1}^N g_i(x) u_i$ for an $N$-player game with necessary assumptions of $f(x)$ and $\{g_i(x)\}_{i=1}^N$. This case and consideration of cooperative games or games with some stochastic element in the dynamics can be directions for the further research.

\bibliographystyle{ieeetr}
\bibliography{ieeepaper}

\vfill

\end{document}